\documentclass[a4paper,11pt]{article}
\usepackage[utf8]{inputenc}
\usepackage{geometry}
\usepackage{amsmath,amssymb}
\usepackage{graphicx,xcolor}
\usepackage{setspace,multicol,tcolorbox}
\usepackage{titlesec}
\usepackage{graphicx} % Required for inserting images
\usepackage{amsmath}
\usepackage{amssymb}
\usepackage{subcaption}
\usepackage{mathtools}
\usepackage{mathrsfs}
\usepackage{bbm}
\usepackage{wrapfig}
\usepackage{lipsum}
\usepackage{caption}
\usepackage{dutchcal}
\usepackage{listings}
\usepackage{soul,color}
\usepackage{enumerate}
\usepackage{lipsum} 
\usepackage[]{mdframed}
\usepackage{tikz}
\usepackage{xparse}
\usepackage{cite}
\usetikzlibrary{shapes ,arrows.meta , positioning}
\DeclareMathAlphabet{\mathcal}{OMS}{cmsy}{m}{n}
\usepackage[numbers, compress]{natbib}  % For numerical citations
\usepackage[breaklinks,colorlinks = true,linkcolor = blue,urlcolor=blue,citecolor=red]{hyperref}

\definecolor{haikya}{RGB}{191, 0, 0}

\titleformat{\section}
  {\bf\Large} % Non-bold, Large font, centered
  {\thesection}                 % Section numbering (I, II, III, ...)
  {1em}                         % Spacing between number and title
  {}

\titleformat{\subsection}
  {\centering \large} % Non-bold, large font, centered
  {\bf\thesubsection}              % Subsection numbering (A, B, C, ...)
  {1em}                         % Spacing between letter and title
  {\bf}
\renewcommand{\thesection}{\Roman{section}}
\newtcolorbox{myboxgreen}{colback=green!15!white,colframe=black!60!green}
\newtcolorbox{myboxred}{colback=red!15!white,colframe=black!50!red}
\newtcolorbox{myboxblue}{colback=blue!20!white}
\newtcolorbox{myboxyel}{colback=yellow!30!white,colframe=yellow!30!white}
\newtcolorbox{myboxbk}{colback=black!10!white,colframe=black!10!white}
\newtcolorbox{mybox}{colback=white,colframe=black}

\title{Levinson's Theorem and its Generalization for Dirichlet L-Functions}
\author{Swapnil Ray\\
\\
supervised by Prof. Soumya Das\\ and co-supervised by Prof. Soumya Bhattacharya\\
}

\date{July 30, 2025}

\begin{document}
\maketitle

\begin{abstract}
In this report, we shall present a proof of Levinson’s theorem, following the ideas due to Matthew P. Young in 2010, which states that one-third of the non-trivial zeroes of the Riemann-zeta function \(\zeta(s)\) lie on the critical line, i.e. the line \(\operatorname{Re}(s) = 1/2\), using a mollified second moment of the zeta-function. Later, we shall present a generalized result for Dirichlet’s L-function presented by Xiaosheng Wu in 2018, using the method of Levinson that more than two-fifths of the non-trivial zeroes of the Dirichlet’s L-functions are on the critical line and in addition, more than two-fifths of the non-trivial zeroes are simple and on the critical line, using a longer mollifier than the one used in the proof of Levinson’s theorem. This is a generalization of a result deduced by Conrey in 1989 that the Riemann zeta-function has at least two-fifth zeroes on the critical line.
\end{abstract}

\tableofcontents

\section{Introduction}
Before starting with the problem which we shall address in this thesis report, we shall discuss a few preliminary concepts.  \\
\subsection{Riemann-zeta function}
        Let $s\in\mathbb{C}$ such that $Re(s)>1$. Then $\zeta(s)$ defined as
     $$\zeta(s):=\sum_{n=1}^{\infty} \frac{1}{n^s}$$ is known as the Riemann-zeta function.
     $\zeta(s)$ can be defined for all $s\in\mathbb{C}$, $Re(s)>0$ and $s\neq1$ as 
         $$\zeta(s)=\frac{s}{s-1}-s\int_{1}^{\infty}(x-\lfloor x \rfloor)x^{-s-1}dx,$$ using the Abel Summation formula.\cite{3}\\
         For various computational purposes, its beneficial to consider the completed zeta function defined as $$\xi(s):=\frac{1}{2}s(s-1)\pi^{-\frac{1}{2}s}\Gamma(\frac{1}{2}s)\zeta(s).$$
         We want to find an analytic continuation of $\xi(s)$ to an entire function of $\mathbb{C}$, and a functional equation which will help us in evaluating $\zeta(s)$ at $\{z:\text{Re(z)}<0\}$.Thus we state the following theorem.\\
         \\
         \textbf{Theorem 1.1:}  $\xi(s)$ satisfies the following statements.\cite{2,1}\\
         \textbf{Functional equation:}\\
         The completed zeta function satisfies the functional equation $$\xi(s)=\xi(1-s)$$.
         \textbf{Analytic continuation:}\\
         For $z\in\mathbb{C}$, with $Re(z)\geq0$ the function $\xi(s)$ can be extended to the entire function $$\xi(s)=\pi^{-s/2}(\frac{s(s-1)}{2})\sum_{n=1}^\infty n^{-s}\Gamma(s/2,\pi n^2z)$$ $$+\pi^{(s-1)/2}(\frac{s(s-1)}{2})\sum_{n=1}^{\infty}n^{s-1}\Gamma((1-2)/2,\pi n^2/z)+\frac{s}{2}z^{(s-1)/2}-\frac{(s-1)}{2}z^{s/2}.$$
         \vspace{0.8 cm}
        \subsection{Zeroes of the zeta function}
        
  \begin{itemize}
      \item $\zeta(s)=0$ for $s=-2,-4,-6,.....$. These are called the trivial zeroes of the Riemann-zeta function. These are basically obtained from the fact that $\Gamma(s)$ has simple poles at negative integers.\\
      \item Besides their are non-trivial zeroes of $\zeta(s)$ which lie in the strip $0\leq Re(s)\leq 1$ called the critical strip.
The zeroes of $\zeta(s)$ are symmetric with respect to the real line.\\
      \item The Riemann hypothesis states that the non-trivial zeroes of the zeta function lie on $Re(s)=1/2,$ known as the critical line. The zeroes lying on the critical line are called critical zeroes of $\zeta(s)$.
      \cite{1}
      \\
   \end{itemize}  
      \subsection{Dirichlet characters}
       For $q\in\mathbb{N}$, $\chi:\mathbb{Z}\rightarrow \mathbb{C}$ is said to be a Dirichlet character modulo $q$ if $\chi$ satisfies\\
(i) $\chi(m)=\chi(n)$ for all $m\equiv n$(mod $q$)\\
(ii)$\chi(mn)=\chi(m)\chi(n)$ for all $m,n\in\mathbb{Z}$.\\
(iii) $|\chi(n)|=1$ for all $n\in\mathbb{Z}$ such that $(n,q)=1$\\
 (iv)$\chi(n)=0$ if and only if $(n,q)>1$\\
 If we define a Dirichlet character $\chi_0$ such that $\chi_o(n)=1$ for all $(n,q)=1$ and $0$ otherwise, then $\chi_0$ is said to be the principle character modulo $q$.\\
 \subsection{Primitive characters}
 Let $\chi$ be a character (mod $q$). We say $d$ is a quasiperiod of $\chi$ if $\chi(m)=\chi(n)$ whenever $m\equiv n$ (mod $d$) and $(mn,q)=1$.
The least quasiperiod of $\chi$ known as the conductor of $\chi$ is a divisor of $q$. \cite{2} \\
    
\textbf{Definition:} A character $\chi$ modulo $q$ is said to be primitive if $q$ is the least quasiperiod of $\chi$.
\subsection{Gauss sums}
For a character $\chi$ modulo $q$, we define the Gauss sum of $\chi$ as
    $$\tau(\chi)=\sum_{a=1}^{q}\chi(a)e(a/q).$$
    \textbf{Properties:}\\
    (i) $\chi(n)\tau(\bar{\chi})=\sum_{a=1}^{q}\bar{\chi}(a)e(an/q)$ for all $(n,q)=1$ and $\bar{\tau(\chi)}=\chi(-1)\tau(\bar{\chi})$.\\
    (ii) $|\tau(\chi)|=\sqrt{q}$ if $\chi$ is a primitive character.\\
    (iii) If $\chi$ is a primitive character modulo $q$, then
        $$\chi(n)=\frac{1}{\tau(\bar\chi)}\sum_{a=1}^{q}\bar{\chi}(a)e(an/q) $$ for any integer $n$.\cite{2}   \\
    \subsection{Dirichlet L-function}
    Let $s\in\mathbb{C}$ with $Re(s)>1$ and $\chi$ be a  character modulo $q$. Then the Dirichlet L-function $L(s,\chi)$ is defined as 
$$L(s,\chi)=\sum_{n=1}^{\infty}\frac{\chi(n)}{n^s}.$$
\textbf{Remark:}\\
For a character $\chi$, let us define $\kappa$ as $\kappa(\chi)=\begin{array}{cc} \{ & \begin{array}{cc}
     0 & \chi(-1)=1  \\
     1 & \chi(-1)=-1
\end{array} \end{array}$\\
Also, we define $\epsilon({\chi})=\frac{\tau(\chi)}{i^{\kappa}\sqrt{q}}$\\
Let $\nu_{\kappa}(z,\chi)=\sum_{-\infty}^{\infty} n^{\kappa}e^{-\pi n^2 z/q}$.\\
    Then, $\nu_{\kappa}(z,\chi)=\frac{\epsilon(\chi)}{z^{1/2+\kappa}}\nu_{\kappa}(1/z,\bar{\chi})$.\\
Let $\chi$ be a primitive character modulo $q$. Then the completed Dirichlet L-function is defined as $\xi(s,\chi)=L(s,\chi)\Gamma((s+\kappa)/2)(q/\pi)^{(s+\kappa)/2}$.\\
\par
The function $\xi(s,\chi)$ can be extended to the entire function
    $$\xi(s,\chi)=(q/\pi)^{(s+\kappa)/2}\sum_{n=1}^{\infty}\chi(n)n^{-s}\Gamma((s+\kappa)/2,\pi n^2 z/q)\\+\epsilon(\chi)(q/\pi)^{(1-s+\kappa)/2}\sum_{n=1}^{\infty}\bar\chi(n)n^{s-1}\Gamma((1-s+\kappa)/2,\pi n^2 /(qz))$$,\\ where $z\in\mathbb{C}$ with $Re(z)\geq 0$.
    Furthermore, $$\xi(s,\chi)=\epsilon(\chi)\xi(1-s,\bar\chi),$$ is the functional equation of $L(s,\chi)$. \cite{2}
\subsection{Zeroes of Dirichlet L-function}
  \begin{itemize}
      \item $L(s,\chi)=0$ for $s=-\kappa,-\kappa-2,-\kappa-4,-\kappa-6,.....$. These are called the trivial zeroes of the Dirichlet L-function.\\
  
      \item The non-trivial zeroes of $\zeta(s)$ lie in the strip $0\leq Re(s)\leq 1$. The strip $Re(s)\in [0,1]$ is called the critical strip.\\

      \item The Generalized Riemann hypothesis states that the non-trivial zeroes of the Dirichlet L-function lie on $Re(s)=1/2.$ \cite{1}
  \end{itemize}
  In the following sections, we shall mention the currently known zero-free region of $\zeta(s)$ from which we shall state the Prime Number Theorem with its error term.
  \\
  \subsection{Zero-free region of $\zeta(s)$}
  There exists a constant $c>0$ such that $\zeta(s)$ has no zeroes in the region $$\sigma\geq 1-\frac{c}{\text{log}(|t|+2)},$$ where $\sigma=\text{Re}(s)$ and $t=\text{Im}(s)$. \cite{4}

   \begin{figure}[h!] 
  \centering
  \includegraphics[width=12cm, height=7cm]{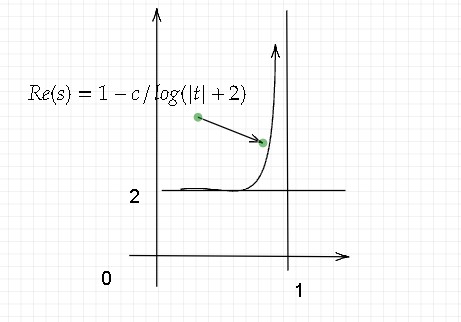}
  \label{fig:Zero free region}
 \end{figure}

 \subsection{Prime Number theorem}
 We are stating the following version of Prime Number theorem without error terms. The proofs are given in \cite{3,4}\\
 \textbf{Theorem 1.2:} Let $\pi(x)$ denote the number of primes $p\leq x$, where $x>0$. Then
     $$\pi(x) \sim \frac{x}{\text{log } x}$$.\\
Let $\psi(x)=\sum_{n\leq x}\Lambda(n)$.
It is proved that the statement of the Prime Number Theorem is equivalent to 
$\psi(x)\sim x$.
We can write $\psi(x)=\frac{1}{2\pi i}\int_{a-i\infty}^{a+i\infty}-\frac{\zeta'(s)}{\zeta(s)}\frac{({{\lfloor{x}\rfloor+\frac{1}{2}}})^s}{s}ds$, where $a>1$ is chosen such that $\zeta(s)\neq 0$ for $Re(s)=a$.\\
Considering the above mentioned zero-free region of $\zeta(s)$ would give us $$\psi(x)=x+O(xe^{-c_{1}(\text{log }x)^{1/2}}),$$ for some positive constant $c_1$.\\
This is the best possible approximation of Prime Number theorem that we have until now.\\
The Riemann hypothesis would give $$\psi(x)=x+O(x^{1/2}\text{log}^2x).$$

\subsection{L-functions}
We shall now discuss the definition and properties of generalized L-functions.\\
\\
\textbf{Definition:} A function $L(s)$ is said to be an L-function if it satisfies the following properties.\\
    (i)$L(s)=\sum_{n\geq1}\lambda(n)n^{-s}=\Pi_{p}(1-\alpha_{1}(p))^{-1}....(1-\alpha_{d}(p))^{-1}$, where $\lambda(1)=1$,$\lambda(n)\in\mathbb{C}$ and the series and the Euler product converge absolutely on $Re(s)>1$.\\
    The $\alpha_{i}(p),1\leq i\leq d$ are called the local roots of $L(s)$ at $p$, and $d$ is called the degree of the Euler product.\\
    (ii) A gamma factor $\gamma(s)=\pi^{-ds/2}\Pi_{j=1}^{d}\Gamma((s+\kappa_{j})/2)$, where the numbers $\kappa_{j}$ are called the local parameters of $L(s)$ at infinity.\\
    (iii) An integer $q\geq1$, called the conductor of $L(s)$, such that $\alpha_{i}(p)\neq 0$ for $(p,q)=1$.
\cite{5}
    \\
    \subsection{Completed L-function}
    The completed L-function is defined as $$\Lambda(s) := q^{
s/2}\gamma(s)L(s).$$     $\Lambda(s)$ admits an analytic continuation to a meromorphic function with at most poles at $s=0$ and $s=1.$ It also satisfies the functional equation $\Lambda(s)=\epsilon\bar\Lambda(1-s)$, where $\bar\Lambda$ is the completion of $\bar{L}(s)=\sum_{n\geq1}\bar{\lambda}(n)n^{-s}$. \\
    \\
    We now state a lemma which will be key to finding an approximate functional equation for $L(s)$.\\
   \textbf{Lemma:} Any L-function $L(s)$ is polynomially bounded in vertical strips $s=\sigma+it$ with $a\leq\sigma\leq b$, $|t|\geq 1$.
   \\
   \\
   We shall now state the approximate functional equation for $L(s)$.\\
   \textbf{Theorem 1.3:} Let $L(s)$ be an L-function. Let $G(u)$ be holomorphic and bounded in $-4<Re(S)<4$, even, and normalized by $G(0)=1$.
    Let $X>0$. Then for $s$ in the strip $0\leq Re(s)\leq 1$ we have
    $$L(s)=\sum_{n}\frac{\lambda(n)}{n^s}V_{s}(\frac{n}{X\sqrt{q}})+\epsilon(s)\sum_{n}\frac{\bar{\lambda}(n)}{n^{1-s}}V_{1-s}(\frac{nX}{\sqrt{q}})+R$$
    where $$V_{s}(y)=\frac{1}{2\pi i}\int_{3-i\infty}^{3+\infty}y^{-u}G(u)\frac{\gamma(s+u)}{\gamma(s)}\frac{du}{u}$$
    and $$\epsilon(s)=\epsilon q \frac{\gamma(1-s)}{\gamma(s)}.$$\\
    If $\Lambda$ is entire then $R=0$, otherwise
    $$R=(res_{u=1-s}+res_{u=-s})\frac{\Lambda(s+u)}{q^{s/2}\gamma(s)}\frac{G(u)}{u}X^{u}.$$
The detailed proofs can be  found in \cite{5}.

\section{Proportion of critical zeroes of $\zeta(s)$}

As of now, there have been various attempts at proving the Riemann hypothesis but it remains an open conjecture till. However, there has been a lot of progress in this direction. In this section, we shall present a short proof of Levinson's theorem proved in 1974, in which he presents a method of computing a lower bound of the proportion of zeroes of $\zeta(s)$ on the critical line, using which he calculates the proportion to be atleast 1/3 of all the non-trivial zeroes.\\
\\
\par
We mention a few notations before proceeding. We let $N(T)$ denote the number of zeroes $s=\sigma+it$ with $0<t<T$, and $N_{0}(T)$ denote the number of such zeroes on the critical line.
We define $\kappa=\text{lim inf}_{T\rightarrow\infty}\frac{N_{0}(T)}{N(T)}$.
The goal is to establish a lower bound on $\kappa$.\\
\subsection{Previously known results}
Before starting with the proof of Levinson's theorem, we shall mention a few key results that were proved before.\\
\begin{itemize}
    \item The first question that comes to our mind is whether or not the number of critical zeroes of $\zeta$ is infinite, as otherwise the proportion would simply be zero. In 1914, it was proved by Hardy that $\zeta(s)$ indeed has infinitely many zeroes on the critical line. Hardy defined the function $$Z(t)=\frac{H(1/2+it)}{|H(1/2+it)|}\zeta(1/2+it),$$ 
    where $H(s)=\frac{1}{2}s(s-1)\pi^{-s/2}\Gamma(s/2).$
    He then showed that $Z$ is real valued and thus, the sign changes of $Z(t)$ correspond to the zeroes of $\zeta(s)$ on the critical line.
    He then proved that $Z(t)$ changes sign infinitely many times by arriving at a contradiction from the fact that finitely many sign changes would imply $\int_{T}^{2T}|Z(t)|dt=\int_{T}^{2T}Z(t)dt$ for large enough T. \cite{8} \\
    \item Although it was established that $\zeta(s)$ has infinitely many zeroes on the critical line, it does not completely eliminate the possibility of the proportion being zero. However it was proved by Selberg in 1942 that the proportion of critical zeroes of $\zeta(s)$ is positive.
    Selberg introduced a mollifier $M(s)$ to use the function $\zeta(s)M(s)$ for detecting zeroes of the zeta function on the critical line. Using this he proved the bound $N_0(T)\gg T\text{log }T$. Using the fact that $N(T)\sim\frac{T}{2\pi}\text{log }T$, he concluded that $\kappa>0$. \cite{8} 
\end{itemize}
    \par
    A method of detecting critical zeroes was introduced by Siegel. He used a formula found in Riemann's notes that
    $$\pi^{-s/2}\Gamma(s/2)\zeta(s)=\pi^{-s/2}\Gamma(s/2)f(s)+\pi^{(1-s)/2}\Gamma((1-s)/2)\bar{f}(1-s),$$
    where $f(s)=\int_{\mathscr{L}}\frac{x^{-s}e^{\pi ix^2}}{e^{\pi ix}-e^{-\pi ix}}dx$. Here $\mathscr{L}$ is a line of slope -1 passing through 1/2. Its easy to see that $\zeta(1/2+it)=0$ whenever the argument of $\pi^{-s/2}\Gamma(s/2)f(s)$ is $\pi/2$ modulo $\pi$. He proved $f(s)$ has $\gg T$ zeroes to the left of the critical line from which he concluded using the argument principle that $N_0(T)\gg T$. \cite{8}\\

\subsection{Levinson method}
Levinson introduced a method which he later used to prove $\kappa> 1/3$. This idea has been used subsequently by other mathematicians while improving the lower bound on $\kappa$. He considered a polynomial $Q(x)$ satisfying $Q(0)=1$ and $V(s)$ defined as
    $$V(s)=Q(-\frac{1}{L}\frac{d}{ds})\zeta(s),$$
    where $L=\text{log}T$. The function $f(s)$ in Siegel's proof was replaced by $V\psi(1-s)$ for a mollifier $\psi$. He used the argument principle to find an upper bound of the zeroes of $V(1-s)$ to the right of the critical line. Levinson concluded that
    $$\kappa \geq 1-\frac{1}{R}\text{log}(\frac{1}{T}\int_1^T|V\psi(\sigma_0+it)|^2dt),$$
    where $\sigma_0=\frac{1}{2}-R/L$ for a positive real number $R$ to be chosen later. The mollifier $\psi$ was defined as
    $$\psi(s)=\sum_{h\leq M}\frac{\mu(h)}{h^{\frac{1}{2}+s-\sigma_{0}}}P(\frac{\text{log }M/h}{\text{log}M}),$$
    where $M=T^{\theta}$ for some $0<\theta<1/2$, $\sigma_0=1/2-R/L$ for some positive real number $R$ and $P(x)$ is a polynomial satisfying $P(0)=0$, $P(1)=1$. The idea was to find an asymptotic formula for a mollified second
    moment of the Riemann-zeta function in order to compute the lower bound on $\kappa$, which we are stating as the following theorem. \\
    \\
    \textbf{Theorem 2.1:}We have $$\frac{1}{T}\int_1^T|V\psi(\sigma_0+it)|^2dt=c(P,Q,R,\theta)+o(1),$$ where
        $$c(P,Q,R.\theta)=1+\frac{1}{\theta}\int_0^1\int_0^1e^{2Rv}(\frac{d}{dx}e^{R\theta x}P(x+u)Q(v+\theta x)|_{x=0})dudv.$$
 As a result of the above mentioned theorem, we have $$\kappa\geq 1-\frac{1}{R}\text{log}(c(P,Q,R,\theta)+o(1))+o(1).$$
    Choosing $P(x)=x$, $Q(x)=1-x$, $R=1.3$, $\theta=0.5$ we compute using Wolfram Mathematica, that $c(P,Q,R,\theta)=2.35...$ and $\kappa\geq 0.35...$\\
Levinson's proof published in 1974 is quite difficult and lengthy. We shall present a brief sketch of a much simpler proof published by Matthew P. Young in 2018. A detailed proof is given in  \cite{6}.\\
\\
\par
Before, we continue, let us choose a smooth function $w(t)$ satisfying\\
   (i) $0\leq w(t)\leq 1$, for all $t\in\mathbb{R}$\\
   (ii) $w$ has compact support in $[T/4,2T]$\\
   (iii) $w^{(j)}(t)=O_{j}(\Delta^{-j})$, for each $j=0,1,2,....,$ where $\Delta=T/L.$

\textbf{Theorem 2.2:} We have
$$\int_{-\infty}^{\infty}w(t)|V\psi(\sigma_{0}+it)|^2dt=c(P,Q,R,\theta)\hat{w}(0)+O(T/L),$$ where
        $$c(P,Q,R.\theta)=1+\frac{1}{\theta}\int_0^1\int_0^1e^{2Rv}(\frac{d}{dx}e^{R\theta x}P(x+u)Q(v+\theta x)|_{x=0})dudv.$$
\\
\\
Theorem 2.1 can be obtained from Theorem 2.2, by choosing $w$ to be an upper bound of the characteristic function of the interval $[T/2,T]$, with support in $[T/2-\Delta,T+\Delta]$.
We get $$\int_{T/2}^T |V\psi(\sigma_0+it)|^2 dt\leq c(P,Q,R,\theta)\hat{w}(0)+O(T/L).$$
In this case $\hat{w}(0)=T/2+O(T/L)$.
Similarly, we get a lower bound and summing over dyadic segments we get the expression.\\
In order to prove Theorem 2.2, we develop a series of lemmas.
\subsection{Completing the proof}
 Let $G(s)=e^{s^2}p(s)$, where $p(s)=\frac{(\alpha+\beta)^2-(2s)^2}{(\alpha+\beta)^2}$, and define
    $$V_{\alpha,\beta}(x,t)=\frac{1}{2\pi i}\int_{1-i\infty}^{1+i\infty}\frac{G(s)}{s}g_{\alpha,\beta}(s,t)x^{-s}ds,$$
    $$g_{\alpha,\beta}(s,t)=\pi^{-s}\frac{\Gamma(\frac{\frac{1}{2}+\alpha+s+it}{2})\Gamma(\frac{\frac{1}{2}+\beta+s-it}{2})}{\Gamma(\frac{\frac{1}{2}+\alpha+it}{2})\Gamma(\frac{\frac{1}{2}+\beta-it}{2})}.$$\\
    Furthermore, set $X_{\alpha,\beta,t}=\pi^{\alpha+\beta}\frac{\Gamma(\frac{\frac{1}{2}-\alpha-it}{2})\Gamma(\frac{\frac{1}{2}-\beta+it}{2})}{\Gamma(\frac{\frac{1}{2}+\alpha+it}{2})\Gamma(\frac{\frac{1}{2}+\beta-it}{2})}$.\\
    \\
    An application of the Theorem 1.3 gives us the following lemma.\\
    \\
    \textbf{Lemma 2.1} If $\text{Re}(\alpha),\text{Re}(\beta)<1/2$, for any $A\geq 0$, we have
            $$\zeta(\frac{1}{2}+\alpha+it)\zeta(\frac{1}{2}+\beta-it)=\sum_{m,n}\frac{1}{m^{\frac{1}{2}+\alpha}n^{\frac{1}{2}+\beta}}(\frac{m}{n})^{-it}V_{\alpha,\beta}(mn,t)$$  $$+X_{\alpha,\beta,t}\sum_{m,n}\frac{1}{m^{\frac{1}{2}-\beta}n^{\frac{1}{2}-\alpha}}(\frac{m}{n})^{-it}V_{-\beta,-\alpha}(mn,t)+O_{A}((1+|t|)^{-A}).$$
    \\
    Using the above mentioned lemma, we arrive at the next lemma, where $w$ is a smooth function satisfying the properties we mentioned as in Theorem 2.2.\\
    \\
    \textbf{Lemma 2.2} Let $h,k$ be positive integers with $h,k\leq T^{2\theta}$, and $\alpha,\beta=O(L^{-1})$. Then
        $$\int_{-\infty}^{\infty}w(t)(\frac{h}{k})^{-it}\zeta(\frac{1}{2}+\alpha+it)\zeta(\frac{1}{2}+\beta-it)dt$$ $$=\sum_{hm=kn}\frac{1}{m^{\frac{1}{2}+\alpha}n^{\frac{1}{2}+\beta}}\int_{-\infty}^{\infty}V_{\alpha,\beta}(mn,t)w(t)dt\\+\sum_{hm=kn}\frac{1}{m^{\frac{1}{2}-\beta}n^{\frac{1}{2}-\alpha}}\int_{-\infty}^{\infty}V_{-\beta,-\alpha}(mn,t)X_{\alpha,\beta,t}w(t)dt$$ $$+O_{A,\theta}(T^{-A})$$
    Let us define $$I_{1}(\alpha,\beta)=\sum_{h,k\leq M}\frac{\mu(h)\mu(k)}{\sqrt{hk}}P(\frac{logM/h}{logM})P(\frac{logM/k}{logM})J_{1},$$
    where $J_{1}$ denotes the first term in the RHS of Lemma 2.2.
    We now state the following lemma.\\
    \\
    \textbf{Lemma 2.3: }For a fixed annuli with $L^{-1}\ll\alpha,\beta\ll L^{-1}$, $|\alpha+\beta|\ll L^{-1}$, we have $$I_{1}(\alpha,\beta)=c_{1}(\alpha,\beta)\hat{w}(0)+O(T/L),$$ where
        $$c_{1}(\alpha,\beta)=\frac{1}{(\alpha+\beta)logM}\frac{d^2}{dxdy}M^{\alpha x+\beta y}\int_{0}^{1}P(x+u)P(y+u)du|_{x=y=0}.$$
    \textbf{Proof:} The Convolution Theorem for Fourier transforms for $1\leq h\leq M$ and\\ $i=1,2,3...$ gives
$$(\frac{logM/h}{logM})^{i}=\frac{i!}{(logM)^{i}}\frac{1}{2\pi i}\int_{1-i\infty}^{1+i\infty}(\frac{M}{h})^{v}\frac{dv}{v^{i+1}}.$$\\
Using the previous two lemmas, we get
$$I_{1}(\alpha,\beta)=\int_{-\infty}^{\infty}w(t)\sum_{i,j}\frac{a_{i}a_{j}i!j!}{(logM)^{i+j}}\sum_{hm=kn}\frac{\mu(h)\mu(k)}{h^{\frac{1}{2}}k^{\frac{1}{2}}m^{\frac{1}{2}+\alpha}n^{\frac{1}{2}+\beta}}$$
$$(\frac{1}{2\pi i})^{3}\int_{1-i\infty}^{1+i\infty}\int_{1-i\infty}^{1+i\infty}\int_{1-i\infty}^{1+i\infty}\frac{M^{u+v}}{h^{v}k^{u}}\frac{g_{\alpha,\beta}(s,t)}{(mn)^{s}}\frac{G(s)}{s}ds\frac{du}{u^{i+1}}\frac{dv}{v^{j+1}}.$$\\
We compute the sum over $h,k,m,n$ as follows
$$\sum_{hm=kn}\frac{\mu(h)\mu(k)}{h^{\frac{1}{2}+v}k^{\frac{1}{2}+u}m^{\frac{1}{2}+\alpha+s}n^{\frac{1}{2}+\beta+s}}=\frac{\zeta(1+u+v)\zeta(1+\alpha+\beta+2s)}{\zeta(1+\alpha+u+s)\zeta(1+\beta+v+s)}A_{\alpha,\beta}(u,v,s),$$
where the arithmetical factor $A_{\alpha,\beta}(u,v,s)$ is given by an absolutely convergent Euler product. Next we move the countours to $Re(u)=Re(v)=\delta$, and then $Re(s)=\delta+\epsilon$, crossing a pole at $s=0$ only since $G(s)$ vanishes at the pole of $\zeta(1+\alpha+\beta+2s)$.
Since $M\leq T^{\theta}$ with $\theta<1/2$ and $t\geq T/2$, the new contour of integration gives $O(T^{1-\epsilon})$ for sufficiently small $\epsilon>0$.\\
    Thus $$I_{1}(\alpha,\beta)=\hat{w}(0)\zeta(1+\alpha+\beta)\sum_{i,j}\frac{a_{i}a_{j}i!j!}{(logM)^{i+j}}J_{\alpha,\beta}(M)+O(T^{1-\epsilon}),$$
    where
    $$J_{\alpha,\beta}(M)=(\frac{1}{2\pi i})^{2}\int_{\epsilon-i\infty}^{\epsilon+i\infty}\int_{\epsilon-i\infty}^{\epsilon+i\infty}M^{u+v}\frac{\zeta(1+u+v)A_{\alpha,\beta}(u,v,0)}{\zeta(1+\alpha+u)\zeta(1+\beta+v)}\frac{du}{u^{i+1}}\frac{dv}{v^{j+1}}.$$
    The proof is completed using the following result, qed.\\
    \\
    \textbf{Result: }  For $\alpha,\beta\ll L^{-1}$,
        $$J_{\alpha,\beta}(M)=\frac{(logM)^{i+j-1}}{i!j!}\frac{d^2}{dxdy}M^{\alpha x+\beta y}\int_{0}^{1}(x+u)^{i}(y+u)^{j}du|_{x=y=0}+O(L^{i+j-2}).$$
    Let us now define $I(\alpha,\beta):=\int_{-\infty}^{\infty}w(t)\zeta(\frac{1}{2}+\alpha+it)\zeta(\frac{1}{2}+\beta-it)|\psi(\sigma_{0}+it)|^2 dt,$ 
where $\psi$ is a mollifier of the form $$\psi(s)=\sum_{h\leq M}\frac{\mu(h)}{h^{\frac{1}{2}+s-\sigma_{0}}}P(\frac{\text{log }M/h}{\text{log }M}),$$
$\sigma_{0}=1/2-R/L$, with $R\ll1$.\\
We now state the following lemma.\\
\\
\textbf{Lemma 2.4: } We have $$I(\alpha,\beta)=c(\alpha,\beta)\hat{w}(0)+O(T/L)$$ for $\alpha,\beta\ll L^{-1}$, where
        $$c(\alpha,\beta)=1+\frac{1}{\theta}\frac{d^2}{dxdy}M^{-\beta x-\alpha y}\int_0^1\int_0^1T^{-v(\alpha+\beta)}P(x+u)P(y+u)du|_{x=y=0}.$$
 \textbf{Proof:} Using the definition of $\psi$ and then applying Lemma 2.2, we write $$I(\alpha,\beta)=I_1(\alpha,\beta)+I_1(-\beta,-\alpha)+O(T^{-A}).$$\\
    We use Lemma 2.3 to complete the proof, qed.\\
\par
Let $Q$ be a polynomial with real coefficients satisfying $Q(0)=1$. Let us define $$V(s)=Q(-\frac{1}{L}\frac{d}{ds})\zeta(s).$$ We are now in a position to complete the proof of Theorem 2.2.\\
    \\
\textbf{Proof of Theorem 2.2:}
    Let us define $I_{smooth}$ to be the LHS.\\
    We observe $I_{smooth}=Q(-\frac{1}{L}\frac{d}{d\alpha})Q(-\frac{1}{L}\frac{d}{d\beta})I(\alpha,\beta)|_{\alpha=\beta=-R/L}$.\\
    Also applying the differential operator to $c(\alpha,\beta)$ evaluated at $\alpha=\beta=-R/L$ gives us $c(P,Q,R,\theta)$, qed.\\
    \\
    \par
    Thus, we have completed the proof of Levinson's Theorem. In the next section, we shall present a generalized result for Dirichlet L-functions. But before that we shall state a few results that have been published after Levinson's theorem, which have increased the lower bound of $\kappa$ .\\
\\
\subsection{Improvements on the lower bound of $\kappa$}
  While computing the lower bound on $\kappa$, we made a few choices regarding the polynomials $P$ and $Q$, and the mollifier $\psi$. While the choice of the polynomials does not change the result significatly, given the initial conditions are satisfied, the choice of the mollifier is a key factor in obtaining the lower bound of $\kappa$. In the mollifier $\psi$ chosen in Levinson's method, there were $\lfloor T^\theta\rfloor$ terms in the summation, where $\theta$ is a chosen positive real number with $\theta<1/2$. However, increasing the value of $\theta$ would result in more number of terms in the mollifier, which will make it easier to detect critical zeroes, and as a result, improve the lower bound of $\kappa$. $\theta=1-\epsilon$ implies the Lindelöf hypothesis which states that for any $\epsilon>0$, $$\zeta(1/2+it)\ll t^\epsilon,$$ and $\theta=\infty$ implies the Riemann hypothesis proved by Bettin and Gonek in 2017\cite{9}. In 1985, Conrey, Heath-Brown and Balasubramaian gave an asymptotic
     formula for I for the Riemann-zeta function which states
     $$I=T\sum_{h,k\leq y}\frac{a(h)\bar{a(k)}}{[h,k]}(\text{log}\frac{T(h,k)^2}{2\pi hk}+2\gamma-1+2\text{log}2)+O(T^{1-\epsilon_\theta}),$$ where $I$ is defined as $$I:=\int_T^{2T} |\zeta(1/2+it)|^2|B(1/2+it)|^2dt,$$
where $\chi$ is a primitive character and $B(s)$ is a Dirichlet polynomial defined as, $$B(s)=\sum_{n\leq y}\frac{a(n)}{n^s}$$ with $a_n \ll n^{\epsilon}$, $y=T^{\theta}$, and $\theta<1$. where $\epsilon_\theta$ was a constant dependent of $\theta$. In general, $\epsilon_\theta>0$ for $\theta<1/2$.
     Besides they were able to increase the value of $\theta$ to $\theta< 9/17$ for the coefficients $a(m)$ in $B(s)$ having the form
     $$a(m)=\mu(m)\mathscr{F}(m) \text{ with } \mathscr{F}\in\mathbb{F}=\{\mathscr{F}:\mathscr{F}(x)\ll_\epsilon x^\epsilon,\mathscr{F'}(x)\ll 1/x \}$$
     which implied the proportion of zeroes on the critical line to be atleast 38\% \cite{10}.
     Later, Conrey extended $\theta$ to 4/7 and proved that the Riemann zeta-function has more than 40.88\% zeroes on the critical line\cite{12}. Later on, Bettin, Chandee and Radziwill extended $\theta$ to $17/33$ \cite{11}. In 2018, Pratt and Robles deduced an asymptotic formula for $\theta<17/33$ \cite{13}. Besides, they extended $\theta$ to $6/11$ for the special coefficient $a(n)=\mu(n)(\mu*\Lambda^{*k})(n)P(\frac{\text{log }y/n}{\text{log y}})$, which improves the least proportion of zeroes of the Riemann zeta-function to be 41.491\%. 

\section{Generalized results for Dirichlet L-functions}
In the previous section, we saw the proof of Levinson's theorem and mentioned some of its subsequent improvisations. We shall now present a similar result for Dirichlet L-functions. In addition, we also obtain a lower bound on the proportion of simple zeroes of a Dirichlet L-function lying on the critical line, The proof we are presenting here was published by Xiaosheng Wu in 2018.\cite{7} Here the upper bound on $\theta$ was taken to be $4/7$ as in Conrey's proof for the case of the Riemann zeta-function in
\cite{12}.
\\
\par
Let $\chi$ be a character modulo $q$ and $L(s,\chi)$ be its associated Dirichlet L-function.
We define $$I(\chi):=\int_T^{2T} |L(1/2+it,\chi)|^2|B(1/2+it,\chi)|^2dt,$$
where $\chi$ is a primitive character and $B(s,\chi)$ is a Dirichlet polynomial defined as, $$B(s,\chi)=\sum_{n\leq y}\frac{\chi(n)a(n)}{n^s}$$ with $a_n \ll n^{\epsilon}$, $y=T^{\theta}$, and $\theta<1$. We denote
 $$\mathscr{L}_{\chi}:=\text{log}\frac{qT}{2\pi}, \text{ and } \mathscr{L}_{1}:=\text{log}\frac{T}{2\pi}.$$
 Let $Q(x)$ be a polynomial. We define
 $$V(s,\chi):=Q(-\frac{1}{\mathscr{L}_{\chi}}\frac{d}{ds})L(s,\chi)$$
 and $$I_R(Q,\chi):=\int_T^{2T} |V(1/2+it+\frac{R}{\mathscr{L}_{\chi}},\chi)|^2|B(1/2+it,\chi)|^2dt$$ for some given real number $R$.\\
Let $N(T,\chi)$ denote the number of zeroes of $L(s,\chi)$ with $0<\sigma<1$ and $|t|\leq T$.
Also let $N_c(T,\chi)$  denote the number of critical zeroes and $N_c^*(T,\chi)$ simple and critical zeroes of $L(1/2+it,\chi)$ with $|t|\leq T$ respectively. We define $$\kappa(\chi)=\frac{N_c(T,\chi)}{N(T,\chi)} \text{ and } \kappa^*(\chi)=\frac{N_c^*(T,\chi)}{N(T,\chi)}.$$
We shall now state our desired result in the following theorem.\\
\\
\textbf{Theorem 3.1:} We have for any Dirichlet character $\chi$,
        $$\kappa(\chi)>0.4172 \text{ and } \kappa*(\chi)> 0.4074$$ for sufficiently large $T$ with $\text{log }q=o(\text{log }T)$.\\
\\
Levinson found a generalized inequality for the proportion of critical zeroes of Dirichlet L-functions which states $$\kappa(\chi)\geq 1-\frac{1}{R}\text{log }(T^{-1}I_R(Q,\chi))+o(1)$$ for a given $R>0$. Also if $Q(x)$ is a linear polynomial, the inequality yields a lower bound for the proportion of simple zeroes on the critical line. To prove Theorem 3.1, we want to find an asymptotic formula for $I_R(Q,\chi)$ quite similar to the case of Riemann-zeta function in order to approximate $I_R(Q,\chi)$ and compute the lower bounds of $\kappa$ and $\kappa^*$ respectively.\\
\par
        Let $\chi$ be a primitive Dirichlet character (mod $q$) with 
 $\text{log}q=o(\text{log}T)$ and $\alpha=a/\mathscr{L}_\chi$, $\beta=b/\mathscr{L}_\chi$ with $a,b\in \mathbb{C}$ and $a,b\ll1$. Let $a(m)\ll_\epsilon m^\epsilon$ for any $\epsilon>0$ with $y=T^\theta$.
We state the following theorem.\\
 \textbf{Theorem 3.2: } We have the following asymptotic relation. $$I_R(Q,\chi)=TQ(\frac{-d}{da})\bar{Q}(\frac{-d}{db}) \{\sum_{h,k\leq y}\frac{\chi_0(hk)(h,k)^{\alpha+\beta}a(h)\bar{a(k)}}{h^{1+\beta}k^{1+\alpha}}$$\\
$$(\frac{2^{1+\alpha+\beta}-1}{1+\alpha+\beta}(\frac{2\pi hk}{qT(h,k)^2})^{\alpha+\beta}L(1-\alpha-\beta,\chi_0)+L(1+\alpha+\beta,\chi_0))\}|_{a=b=-R}+O(R^{1-\epsilon_0}).$$  
Here $\epsilon_\theta$ is a constant depending on $\theta$ as follows:\\
    (A) We have $\epsilon_\theta>0$ for all $\theta<17/33$;\\
    (B) We have $\epsilon_\theta>0$ for all $\theta<4/7$ when $a(n)=\mu(n)(\mathscr{F}_0+\mathscr{F}_1.(\mathscr{F}_2*\mathscr{F}_3))(n)$ with $\mathscr{F}_i$ separable in $\mathbb{F}=\{\mathscr{F}:\mathscr{F}(x)<<_\epsilon x^\epsilon, \mathscr{F'}(x)= <<\frac{1}{x} \}$ for $0\leq i\leq 3.$ In addition, it also holds when one of $\mathscr{F}_2$ and $\mathscr{F}_3$ is separable in $\mathbb{F}=\{\mathscr{F}:\mathscr{F}(x)<<_\epsilon x^\epsilon, \mathscr{F}(x)=0 \text{ for } x>y^{3/4} \}$ and other $\mathscr{F}_i$ are separable in $\mathbb{F}$. \\
    \\
    \textbf{Definition: } Let $\mathbb{S}$ be set of arithmetic functions. We say $\mathscr{F}$ is separable in $\mathbb{S}$ if $\mathscr{F}\in\mathbb{S}$ and $\mathscr{F}(mn)$ can be expressed as a finite sum of $\mathscr{F}_i(m)\mathscr{F}_j(n)$ with $\mathscr{F}_i,\mathscr{F}_j\in\mathbb{S}$.
    \\
    \par
    
    We consider $\eta>0$, $\Delta=T^{1-\eta}$, $y=T^{\theta}$ and $s_0=1/2+iw$ with $T\leq w\leq 2T$.
    Let us define $$g(\alpha,\beta,w):=\frac{1}{i\Delta\pi^{1/2}}\int_{(1/2)}e^{(s-s_0)^2\Delta^{-2}}L(s+\alpha,\chi)L(1-s+\beta,\bar{\chi})$$\\ $$B(s,\chi)B_1(1-s,\bar{\chi})ds,$$
    where $B_1(1-s,\bar{\chi})=\sum_{n\leq y}\frac{\bar{\chi(n)}\bar{a(n)}}{n^{1-s}}$. We state the following proposition from which Theorem 3.2 follows.
    \\
    \\
    \textbf{Proposition: } $g(\alpha,\beta,w)$ satisfies the following asymptotic relation $$g(\alpha,\beta,w)=\sum_{h,k\leq y}\frac{\chi_0(hk)(h,k)^{\alpha+\beta}a(h)\bar{a(k)}}{h^{1+\beta}k^{1+\alpha}}$$\\
        $$(\frac{2\pi hk}{qT(h,k)^2})^{\alpha+\beta}L(1-\alpha-\beta,\chi_0)+L(1+\alpha+\beta,\chi_0))+O(T^{-\epsilon_0}).$$
    (A1) $\epsilon_\theta>0$ for any given $\theta<17/33$.\\
    (B1) When $a(n)$ has a special form as in (B) of Theorem 3.2, $\epsilon_\theta>0$ for any given $\theta<4/7$.\\
    \\
    Next, we shall present a proof of the proposition which completes the proof of Theorem 3.2.\\
    \subsection{Main term of the proposition}
    The main term of the proposition is obtained from some pivotal lemmas, which we shall state below. The following lemma is obtained as a result of the Convolution theorem for Fourier transforms.\\
    \\
    \textbf{Lemma 3.1: } Let $1<c<2$. Then 
        $$\frac{1}{i\Delta\pi^{1/2}}\int_{(c)}e^{(s-z)^2\Delta^{-2}}\Gamma(s)(2\pi ix)^{-s}ds=\int_0^\infty v^z\text{exp }(-\frac{\Delta^2\text{log }^2v}{4})e(-xv)\frac{dv}{v},$$ for any $x\neq0,z$ amd $\Delta>0$.\\
    \\
    \par
    Now, let us define $$J(x,s_0,\beta,\Delta,\chi)=\frac{1}{i\Delta\pi^{1/2}}\int_{(c)}e^{(s-s_0)^2\Delta^{-2}}H(1-s+\beta,\chi)x^{-s}ds,$$ where $1<c<2$ and $q,\Delta,s_0,\beta$ are defined as in the proposition, and $$H(1-s,\chi)=(2\pi)^{-s}q^{s-1}\tau(\chi)\Gamma(s)(e^{-\pi is/2}+\chi(-1)e^{\pi is/2}).$$
    The following lemma is obtained as an application of Lemma 3.1.\\
    \\
    \textbf{Lemma 3.2: } We have $$J=\frac{\tau(\chi)}{qx^{\beta}}\int_0^\infty v^{s_0-\beta}\text{exp }(-\frac{\Delta^2\text{log }^2v}{4})$$\\
        $$(e(-\frac{xv}{q})+\chi(-1)e(\frac{xv}{q}))\frac{dv}{v},$$ for any $x\neq 0$.\\
    \\
    
    Let us consider integers $H,K$ with $K\geq 1$ such that any two of $q,H,K$ are co-prime.
    Suppose $\alpha,\beta,s \in\mathbb{C}$. We define
    $$D(s,\alpha,\beta,\frac{H}{Kq},\chi):=\sum_{m,n}\frac{\chi(m)\chi(n)}{m^{s+\alpha}n^{s+\beta}}e(\frac{mnH}{Kq}).$$
    We obtain an analytic continuation and the functional equation of $D$ from the next lemma which we shall state below. \\
    \\
    \textbf{Lemma 3.3: } The function $$D(s,\alpha,\beta,\frac{H}{Kq},\chi)-K^{1-2s-\alpha\beta}\tau(\chi)\chi(K)\bar{\chi}(H)$$
        $$(q^{s-\alpha}L(s+\beta,\chi_0)\zeta(s+\alpha)+q^{-s-\beta}L(s+\alpha,\chi_0)\zeta(s+\beta)$$ $$-q^{-2s-\alpha-\beta}\phi(q)\zeta(s+\alpha)\zeta(s+\beta))$$
        is an entire function of $s$. Also $D$ satisfies
        $$D(1-s,\alpha,\beta,\frac{H}{Kq},\chi)=\frac{2}{(Kq)^2}(\frac{Kq}{2\pi})^{2s-\alpha-\beta}\Gamma(s-\alpha)\Gamma(s-\beta)$$
        $$\{ \text{cos }\frac{\pi}{2}(2s-\alpha-\beta)A_1(s,\alpha,\beta,\frac{H}{Kq},\chi)+\{ \text{cos }\frac{\pi}{2}(\alpha-\beta)A_2(s,\alpha,\beta,\frac{H}{Kq},\chi)\},$$
        with $$A_1(s,\alpha,\beta,\frac{H}{Kq},\chi)=\sum_{1\leq v,u \leq Kq}\chi(u)\chi(v)e(\frac{uvH}{Kq})F(s-\alpha,\frac{u}{Kq})F(s-\beta,\frac{v}{Kq}),$$ and 
        $$A_2(s,\alpha,\beta,\frac{H}{Kq},\chi)=\sum_{1\leq v,u \leq Kq}\chi(u)\chi(v)e(\frac{uvH}{Kq})F(s-\alpha,\frac{u}{Kq})F(s-\beta,\frac{-v}{Kq}),$$
        where $F(s,\chi)=\sum_{n=1}^\infty e(nx)n^{-s}$.
        Moreover, $D(0,\alpha,\beta,\frac{H}{Kq},\chi)\ll_\epsilon q^{3/2+\epsilon}K^{1+\epsilon}$ for any $\epsilon>0$ when $\alpha,\beta\ll(\text{log }Kq)^{-1}$.
    Let us now assume $\alpha,\beta,x\in\mathbb{C}$ with $\alpha\neq\beta,\alpha,\beta\neq1, Im(x)>0$.
Let us define $$S(x,\alpha,\beta,\frac{H}{Kq},\chi):=\sum_{m,n}\frac{\chi(m)\chi(n)}{m^\alpha n^\beta}e(\frac{mnH}{Kq})e(mnx),$$ where $q,H,K$ are chosen as in Lemma 3.3.
  Using the Convolution theorem for Fourier transforms and Lemma 3.3, we obtain the following lemma.\\
  \\
  \textbf{Lemma 3.4: }For any $c>1-\text{min}\{Re(\alpha),Re(\beta)\}$,
    $$S(x,\alpha,\beta,\frac{H}{Kq},\chi)$$
    $$=L(1-\alpha+\beta,\chi_0)K^{-1+\alpha-\beta}q^{-1}\tau(\chi)\bar{\chi}(H)\chi(K)z^{-1+\alpha}\Gamma(1-\alpha)$$
    $$+L(1+\alpha-\beta,\chi_0)K^{-1-\alpha+\beta}q^{-1}\tau(\chi)\bar{\chi}(H)\chi(K)z^{-1+\alpha}\Gamma(1-\beta)$$
    $$+D(0,\alpha,\beta,\frac{H}{Kq},\chi)+\frac{1}{(Kq)^2 \pi i}\int_{(c)}z^{s-1}\Gamma(1-s)\Gamma(s-\alpha)\Gamma(s-\beta)(\frac{Kq}{2\pi})^{2s-\alpha-\beta}$$
    $$\{ \text{cos }\frac{\pi}{2}(2s-\alpha-\beta)A_1(s,\alpha,\beta,\frac{H}{Kq},\chi)+\{ \text{cos }\frac{\pi}{2}(\alpha-\beta)A_2(s,\alpha,\beta,\frac{H}{Kq},\chi)\}ds$$
    with $z=2\pi ix$.\\
    \\
    \par
    Let $\alpha,\beta\ll\text{log }^{-1}T$, $0<\delta<\pi/2$, $z=1/2+\beta+iw$ with $T\leq \omega\leq 2T$ and $\Delta=T^{1-\eta}$ with $\eta>0$.
We define $$r_{\delta}(z,\alpha):=\int_{L_{\delta}}v^z\text{exp}(-\frac{\Delta^2\text{log}^2v}{4})(v-1)^{-1+\alpha}\frac{dv}{v},$$
where $L_\delta$ is the half-line $L_\delta=\{re^{i\delta}:r>0\}$.\\
\par
 Let $$W(z,\alpha)=\Gamma(1-\alpha)\{(-2\pi i)^\alpha r_\delta(z,\alpha)-(2\pi i)^{\alpha}r_{-\delta}(z,\alpha)\}.$$ 
We now state one more lemma that we need to obtain the main term of the proposition. This lemma is obtained as an application of Cauchy's theorem and Lemma 3.1.\\
\\
\textbf{Lemma 3.5: } $W(z,\alpha)$ satisfies the following asymptotic relation.
$$W(z,\alpha)=-2\pi i(\frac{w}{2\pi})^{-\alpha}+O(T^{-\eta}).$$
\textbf{Completing of the proof: }
We use the functional equation for $L(1-s+\beta,\chi)$ and Lemma 3.2 to interchange the summation and integration and write 
    $$g(\alpha,\beta,w)=\frac{\tau(\bar{\chi})}{q}\sum_{h,k\leq y}\frac{\chi_0(hk)(h,k)^{\alpha+\beta}a(h)\bar{a(k)}}{h^{1+\beta}k^{1+\alpha}}\sum_{m,n}\chi(m)\chi(n)m^{-\alpha}n^\beta$$
        $$\int_0^\infty v^{s_0-\beta}\text{exp }(-\frac{\Delta^2\text{log }^2v}{4})
        (e(-\frac{mnhv}{q})+\bar{\chi}(-1)e(\frac{mnhv}{q}))\frac{dv}{v}$$
        $$+O(\text{exp}(-T^{2\eta})).$$
\par
Let $\delta>0$ and $L_\delta$ be the half-line as in Lemma 3.5. We express the above integral as a sum of two integrals and use Cauchy's theorem to move one path to $L_\delta$ and one to $L_{-\delta}$ to have
    $$g(\alpha,\beta,w)=\frac{\tau(\bar{\chi})}{q}\sum_{h,k\leq y}\frac{\chi_0(hk)(h,k)^{\alpha+\beta}a(h)\bar{a(k)}}{h^{1+\beta}k^{1+\alpha}}(\bar{\chi}(-1)I_1+I_2)+O(\text{exp}(-T^{2\eta})).$$
    We now use Lemma 3.4 to split both $I_1,I_2$ as $$I_i=M_i+R_i+E_i$$ for $i=1,2$.
    We apply Lemma 3.5 to $M_1,M_2$ to obtain the main term of the proposition. We shall now present a brief description of the error term of the proposition.\\
    \subsection{Error term of the proposition}
     The error term is obtained from the terms containing $R_i$ and $M_i$, with $i=1,2$ in the asymptotic formula of $g(\alpha,\beta,w)$ as in the proof of the main term of the proposition. We use the upper bound for $D(0,\alpha,\beta,\frac{H}{Kq},\chi)$ given in Lemma 3.3 to obtain $R_i\ll T^{-20}$. After some long calculation using various other previously known bounds, we prove $$Z\ll_\epsilon y^{7/8}T^{-1/2+11/2\eta+\epsilon}+y^{7/4}T^{-1+11/2\eta+\epsilon},$$ where $Z$ is the term in the expression of $I_i$ involving $E_i$. This completes the proof of the proposition. The proof of the error term is quite lengthy and is given in more details in \cite{7}.\\

    \subsection{Proof of Theorem 3.1} We shall now use Theorem 3.2 to obtain a lower bound on $\kappa(\chi)$ and $\kappa^*(\chi)$. As mentioned earlier, Levinson method for Dirichlet L-functions yields the inequality $$\kappa(\chi)\geq 1-\frac{1}{R}\text{log }(T^{-1}I_R(Q,\chi))+o(1)$$ for a suitable $R>0$.
    If $Q(x)$ is a linear polynomial, the above inequality yields a lower bound for $\kappa^*(\chi)$.\\
    \par
    To prove Theorem 3.1, we consider $$B(s,\chi)=\sum_{n\leq y}\frac{\chi(n)a(n)}{n^{s}}$$ where $y=T^\theta$ with $\theta=4/7-\epsilon$ and the coefficients
    $$a(n)=\mu(n)(P_1(\frac{\text{log }y/n}{\text{log }y})+P_2(\frac{\text{log }y/n}{\text{log }y})\sum_{p|n,p\leq y^{3/4}}P(\frac{\text{log }y/n}{\text{log }y})),$$ where $P_1,P_2,P$ are real polynomials with $P_1(0)=P_2(0)=P(0)=0$ and $P_1(1)=1$.  Since $\alpha,\beta\ll\mathscr{L}_\chi^{-1}$, $(\frac{2\pi}{qT})^{\alpha+\beta}=e^{-a-b}$, it follows from Theorem 3.2 that $$I_R(Q,\chi)\sim TQ(\frac{-d}{da})\bar{Q}(\frac{-d}{db})\frac{\Sigma(\beta,\alpha)-e^{-a-b}\Sigma(-\alpha,-\beta)}{\alpha+\beta}|_{a=b=-R},$$ with 
    $$\Sigma(\alpha,\beta):=\sum_{h,k\leq y}\frac{\chi_0(hk)(h,k)^{\alpha+\beta}a(h){a(k)}(h,k)^{1+\alpha+\beta}}{h^{1+\beta}k^{1+\alpha}}.$$
    After some calculations using the lemmas we obtained in this section and some previously known bounds due to Mertens Theorem and Levinson (1974), we arrive at the following approximation $$\Sigma(\alpha,\beta)=\frac{q}{\phi(q)\theta\mathscr{L}}\mathscr{M}+O((\text{log log})^7\text{log}^{-2}y),$$ where $\mathscr{M}$ is a positive real number obtained based on the choice of $P_1,P_2,P$.\\
    \\
    \textbf{Proportion of critical zeroes:}\\
    We choose $\theta=4/7-\epsilon$, $R=1.3$ and 
   $$Q(x)=1-0.642x-1.227(x^2/2-x^3/3)-5.178(x^3/3-x^4/2+x^5/5),$$
   $$P_1(x)=x-0.617x(1-x)-0.125x^2(1-x)-0.148x^3(1-x),$$
   $$P_2(x)=x,$$
   $$P(x)=1.55x-1.564x^2+0.177x^3,$$ to obtain $K(\chi)>0.4172$.\\
   \par
   We take $R=1.116$ and $$Q(x)=1-1.032x,$$
   $$P_1(x)=x-0.525x(1-x)-0.183x^2(1-x)-0.085x^3(1-x),$$
   $$P_2(x)=x,$$
   $$P(x)=0.838x-0.938x^2-0.084x^3,$$ to obtain $\kappa^*(\chi)>0.4074$.\\
   Thus we prove Theorem 3.1.\\

   \section{Conclusion} In this report, we saw a simplified proof of Levinson's theorem by Matthew P. Young, and a generalization of Conrey's result to Dirichlet L-functions with the upper bound of $\theta$ to be $4/7$. However, I feel like it is a responsibility on my behalf to mention that further results have been deduced for Dirichlet L-functions. Conrey, Soundararajan and Iwaniec proved in 2013 that the lower bound of $\kappa(\chi)$ for Dirichlet L-functions is 0.5865.\cite{14} Furthermore, it was proved by Xiaosheng Wu that the proportion of simple and critical zeroes is atleast 0.60261 for Dirichlet L-functions  in 2016.\cite{15}

   \bibliography{reference}

\end{document}